\documentclass[a4paper,fleqn]{cas-sc}

\usepackage[authoryear,longnamesfirst]{natbib}

\usepackage{framed}
\usepackage{subfigure}

\newtheorem{lemma}{\bf Lemma}
\newtheorem{proposition}{\bf Proposition}
\newtheorem{theorem}{\bf Theorem}
\newtheorem{assumption}{\bf Assumption}
\newtheorem{remark}{\bf Remark}

\def\tsc#1{\csdef{#1}{\textsc{\lowercase{#1}}\xspace}}
\tsc{WGM}
\tsc{QE}
\tsc{EP}
\tsc{PMS}
\tsc{BEC}
\tsc{DE}

\begin{document}
\let\WriteBookmarks\relax
\def\floatpagepagefraction{1}
\def\textpagefraction{.001}
\shorttitle{Research on the descent direction of PC algorithms for pseudo-convex/convex optimization problems}
\shortauthors{Ting Li, Deren Han, Tanxing Wang, Xingju Cai}

\title [mode = title]{Research on the descent direction of prediction correction algorithms for pseudo-convex/convex optimization problems}
\tnotemark[1]

\tnotetext[1]{This research is supported by the National Natural Science Foundation of China (Grant Nos.12131004, 12471290, 12126603), Grant 2021YFA1003600, Scientific and Technological Innovation Project of China Academy of Chinese Medical Sciences, Grant No.CI2023C064YLL.}


\author[1]{Ting {Li}}[
                        orcid=0009-0009-4142-3492]

\ead{liting_math@buaa.edu.cn}


\affiliation[1]{organization={School of Mathematical Sciences,LMIB},
                addressline={Beihang University},
                city={Beijing},
                postcode={100191},
                country={China}}

\author[1]{Deren {Han}}[]
\ead{handr@buaa.edu.cn}

\author[2]{Tanxing {Wang}}[%
]
\ead{wangtanxing2021@163.com}
\author[2]{Xingju {Cai}}[%
   ]
\cormark[1]
\ead{caixingju@njnu.edu.cn}

\credit{Data curation, Writing - Original draft preparation}

\affiliation[2]{organization={School of Mathematical Sciences, Key Laboratory for NSLSCS, Ministry of Education, Nanjing Normal University},
                city={Nanjing},
                postcode={210023},
                postcodesep={},
                country={China}}

%

\cortext[cor1]{Corresponding author}
%

\begin{abstract}
Prediction-correction algorithms are a highly effective class of methods for solving pseudo-convex optimization problems. The descent direction of these algorithms can be viewed as an adjustment to the gradient direction based on the prediction step. This paper investigates the adjustment coefficients of  these descent directions and offers explanations from the perspective of differential equations. Unlike existing algorithms where the adjustment coefficient is always set to 1, we establish that the range of the adjustment coefficient lies within (1/2,1] for pseudo-convex optimization problems, and [0,1] for convex optimization problems. We also provide rigorous convergence proofs for these proposed algorithms. Numerical experiment results show that the algorithms perform best when the value of the adjustment coefficient makes the algorithm approach or equal to those in differential equations with higher-order global discrete error.
\end{abstract}

%

\begin{keywords}
Prediction correction algorithm \sep Judiciously adjustment coefficient \sep Adaptive stepsize \sep Discrete differential equations \sep Trapezoidal formula
\end{keywords}

\maketitle

\section{Introduction}

This paper first focuses on the study of pseudo-convex optimization problems, which are more general than convex optimization problems. For instance, fractional programming problems, commonly encountered in fields such as wireless communications, economics, management, information theory, optics, graph theory, and computer science \cite{BM19,CY13,MT20,SM92,ZJ15}, are often pseudo-convex rather than convex optimization problems \cite{Ma70}. Directly developing numerical algorithms for pseudo-convex optimization problems can circumvent the complexities and potential errors associated with transforming the model into a convex problem \cite{SY18}.

Prediction-correction algorithms are a highly effective class of methods for solving pseudo-convex optimization problems \cite{CH13,CL17,FP03,HY04,HY13}. In this paper, we express the descent direction of prediction-correction algorithms as an adjustment of the gradient direction using information from the prediction step. We investigate the range and the optimal value of the adjustment coefficient. In many existing well-performing algorithms,  the adjustment coefficient is typically fixed at 1, such as the extra-gradient algorithm \cite{JL18,Ko76}, the improved extra-gradient algorithm \cite{HL02}, the projection contraction algorithm \cite{HY04}, and the forward-backward algorithm \cite{Ts00}. This paper addresses the question: if this coefficient changes, what would be the range of the coefficient and which one is the most efficient? We provide an answer and offer the convergence analysis. We highlight that changing the adjustment coefficient effectively rotates the descent direction and provide an illustrative figure for better understanding. Additionally, considering that adaptive step sizes can significantly improve the convergence speed of the algorithm \cite{HL02,DC18,OA22}, we also employ the line search strategies for the step sizes in the proposed algorithms.

In recent years, there has been growing interest in the relationship between differential equations and optimization problems. Scholars have already mapped many optimization algorithms to different discrete forms of differential equations \cite{AC20,CL21,DO17,PB13,SB16,ZB20}. To provide a deeper theoretical understanding of the proposed algorithms, we explain several specific adjustment coefficients from the perspective of discrete differential equations.


Some pseudo-convex models can be conveniently and efficiently transformed into convex models using techniques such as the Charnes-Cooper transform and Dinkelbach’s transform in fractional programming problems \cite{SY18}. For the resulting convex optimization problems, we design a new prediction-correction algorithm where the correction step can have a different step size from the prediction step, provided that the adjustment coefficient can also be modified. We explain that several specific adjustment coefficients in this algorithm are also closely related to methods in discrete differential equations.


The main contributions of this paper are as follows:
\begin{itemize}
	\item By considering the descent direction of the prediction-correction algorithm as an adjustment to the gradient direction, we propose new prediction-correction algorithms that allow for altering the adjustment coefficient. We illustrate through a figure that modifying the adjustment coefficient essentially rotates the descent direction. From the perspective of discrete differential equations, we explain that when the adjustment coefficient is $\frac{1}{2}$, the algorithm corresponds to the trapezoidal formula, which has a higher-order global discrete error. This also provides an explanation for the optimal adjustment coefficient.
	\item We establish that the range of the adjustment coefficient lies within (1/2,1] for pseudo-convex optimization problems, and [0,1] for convex optimization problems. We also provide rigorous convergence proofs for these proposed algorithms.
	\item The numerical experiment demonstrates the effectiveness of our algorithms and shows that it is better to make the adjustment coefficient close to $\frac{1}{2}$ while ensuring convergence. This is also consistent with the explanation from the perspective of differential equations.
\end{itemize}

The rest of this article is arranged as follows. Section 2 is divided into two subsections. The first subsection is some basic notations and properties. The second subsection provides the motivation for the proposed algorithms and offers an explanation from the perspective of discrete differential equations for some specific adjustment coefficients. Section 3 is divided into three subsections. Subsections 3.1 and 3.2 respectively introduce the improved prediction-correction algorithms for solving pseudo-convex and convex optimization problems, along with their convergence analyses. Subsection 3.3 explains that, under certain conditions, the algorithm presented in subsection 3.2 can be fully aligned with the algorithm derived from the perspective of discrete differential equations. The numerical experiments are given to show the superiority of the algorithms proposed in this paper in Section 4. Section 5 is the conclusion.

\section{Preliminaries}
\subsection{Basic notions and properties}

The optimization problem considered in this paper is
\begin{eqnarray}\label{eqn1-2}
	\min_{x\in\mathbb{R}^n}f(x),
\end{eqnarray}
where $f: \mathbb{R}^{n}\rightarrow \mathbb{R}$ is a smooth pseudo-convex or convex function.

We next give some basic notions, definitions, and properties that are useful in the subsequent paper.

For vectors $x,y\in\mathbb{R}^{n}$, $\langle x,y\rangle$ is the standard inner product. $\|x\|=\sqrt{\langle x,x\rangle}$ is the Euclidean norm. We use $\textbf{P}_{[a,b]}(x)$ to denote the projection of $x$ on the interval $[a,b]$. A differentiable function $f: \mathbb{R}^{n}\rightarrow \mathbb{R}$ is said to be convex, if it satisfies
\begin{eqnarray*}\label{eqn2-1}
	f(y)\geq f(x)+\langle\nabla f(x),y-x\rangle, \forall x,y\in \mathbb{R}^{n}.
\end{eqnarray*}
An operator $T: \mathbb{R}^{n}\rightarrow \mathbb{R}^n$ is said to be monotone, if for all $x, y \in \mathbb{R}^{n}, \langle T(y)-T(x),y-x\rangle\geq 0$. An operator $T: \mathbb{R}^{n}\rightarrow \mathbb{R}^n$ is said to be pseudo-monotone, if for all $x, y \in \mathbb{R}^{n}, \langle T(x),y-x\rangle\geq 0$ implies $\langle T(y),y-x\rangle\geq 0$.  A differentiable function $f: \mathbb{R}^{n}\rightarrow \mathbb{R}$ is said to be pseudo-convex if its gradient is pseudo-monotone. An operator $T: \mathbb{R}^{n}\rightarrow \mathbb{R}^n$ is said to be $L$-Lipschitz continuous, if
\begin{eqnarray*}\label{eqn2-2}
	\|T(y)-T(x)\|\leq L\|y-x\|, \forall x,y\in \mathbb{R}^{n}.
\end{eqnarray*}
Let $x^*$ be a stationary point of the unconstrained differentiable optimization problem (\ref{eqn1-2}), i.e., $\|\nabla f(x^*)\|=0$, and define the stationary point set $X^*=\{x^*:\|\nabla f(x^*)\|=0\}$. The following lemma gives the property of the gradient of a differentiable convex function.
\begin{lemma}\label{lem2-1}(\cite{Ne18}, Theorem 2.1.3)
	For a differentiable convex function $f: \mathbb{R}^{n}\rightarrow \mathbb{R}$, its gradient operator $\nabla f$ is monotone, and thus pseudo-monotone.
\end{lemma}
The following descent lemma is well-known and useful for our analysis of convergence later.
\begin{lemma}\label{lem2-2}(\cite{Ne18}, Theorem 2.1.5)
	If $\nabla f: \mathbb{R}^{n}\rightarrow \mathbb{R}^n$ is $L$-Lipschitz continuous, then for all $x, y\in \mathbb{R}^{n}$,
	\begin{eqnarray*}
		f(x)\leq f(y)+\langle\nabla f(y),x-y\rangle+\frac{L}{2}\|x-y\|^2.
	\end{eqnarray*}
\end{lemma}
Under the convexity of $f$ and the Lipschitz continuity of $\nabla f$, an inequality can be derived in the following Lemma.
\begin{lemma}\label{lem2-3}(\cite{Ne18}, Theorem 2.1.5)
	If $f$ is a convex function and $\nabla f: \mathbb{R}^{n}\rightarrow \mathbb{R}^n$ is a $L$-Lipschitz continuous operator, then for all $x, y\in \mathbb{R}^{n}$,
	\begin{eqnarray*}
		\langle\nabla f(x)-\nabla f(y),x-y\rangle\geq\frac{1}{L}\|\nabla f(x)-\nabla f(y)\|^2.
	\end{eqnarray*}
\end{lemma}

\subsection{Motivation}

The prediction-correction algorithms is a highly effective method for solving pseudo-convex optimization problems. We analyse some typical prediction-correction algorithms which use the gradient step: $z^k=x^k-h_k\nabla f(x^{k})$ as a prediction to illustrate the motivation behind the proposed algorithms.

For unconstrained optimization problem (\ref{eqn1-2}), the correction step of the extra-gradient algorithm \cite{JL18,Ko76} and forward-backward-forward algorithm \cite{Ts00} can be both written as
\begin{eqnarray}\label{eqn1-71}
	x^{k+1}_{EG/FBF}=x^k-h_k\left(\nabla f(x^{k})-(\nabla f(x^{k})-\nabla f(z^{k}))\right).
\end{eqnarray}
The numerical performance of the improved extra-gradient algorithm \cite{HL02} and the projection contraction algorithm \cite{HY04} remains difficult to surpass. Their correction steps can be expressed as:
\begin{eqnarray}\label{eqn1-18}
	x^{k+1}_{IEG/PC}=x^k-\eta\alpha_kh_k\left(\nabla f(x^{k})-(\nabla f(x^{k})-\nabla f(z^{k}))\right).
\end{eqnarray}
We find that algorithms (\ref{eqn1-71}), (\ref{eqn1-18}) both use direction $\nabla f(x^{k})-(\nabla f(x^{k})-\nabla f(z^{k}))$ as the descent direction, where $\nabla f(x^{k})-\nabla f(z^{k})$ is used to adjust the gradient direction $\nabla f(x^{k})$, and the adjustment coefficient is 1. A natural question is what is the best adjustment coefficient? If we reset $\beta(\nabla f(x^{k})-\nabla f(z^{k}))$ as the adjustment direction, what is the range and the best choice of $\beta$? Thus, this paper considers the following improved prediction correction algorithm whose correction step is
\begin{eqnarray}\label{eqn1-13}
	x^{k+1}_{IPC}=x^k-h_k\left(\nabla f(x^{k})-\beta(\nabla f(x^{k})-\nabla f(z^{k}))\right).
\end{eqnarray}
We provide the range of $\beta$ and give the convergence analysis. For the optimal value of $\beta$, we provide an explanation from the perspective of differential equations. Regarding the step size $h_k$, we give two selection ways: constant step size and adaptive step size. For the convex optimization problem, we propose the following form of the algorithm with correction step written as
\begin{eqnarray}\label{eqn1-20}
	x^{k+1}=x^k-\eta\alpha_kh_k\left(\nabla f(x^{k})-\beta(\nabla f(x^{k})-\nabla f(z^{k}))\right).
\end{eqnarray}

\begin{figure}[htbp]
	\begin{center}
	\includegraphics[width=4.5in]{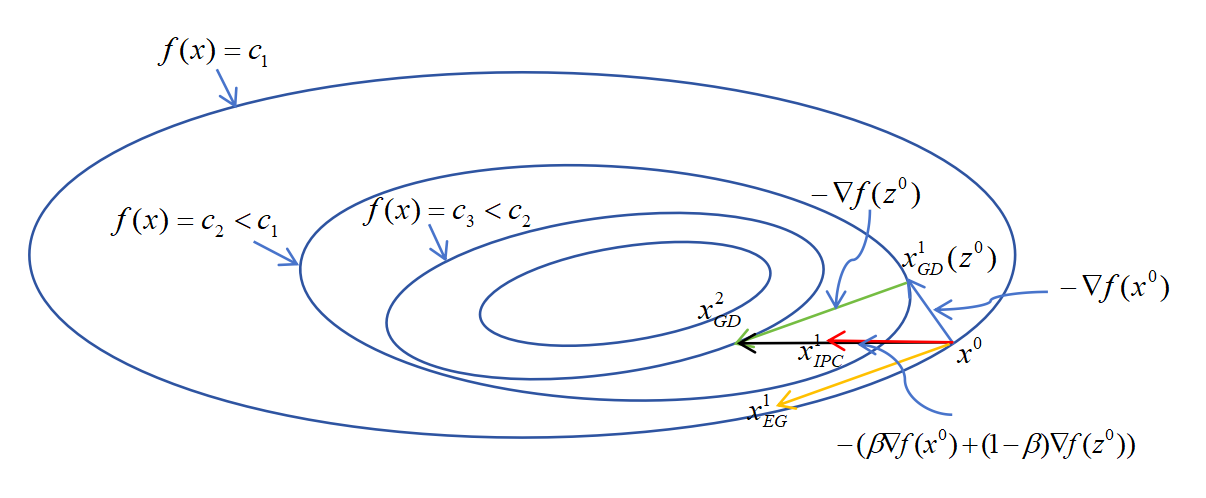}
\caption{Comparison of iterating one step with different adjustment coefficient in Algorithm (\ref{eqn1-13})}
	\label{fig3}		
	\end{center}

\end{figure}

The difference of iterating one step with different adjustment coefficient $\beta$ in Algorithm (\ref{eqn1-13}) can be intuitively shown in Figure \ref{fig3}.
In this figure, $x_{GD}^{1}$ corresponds to $\beta =0$; $x_{EG}^{1}$ corresponds to $\beta =1$; $x_{IPC}^{1}$ corresponds to $\beta \in(0,1)$. As can be seen from this figure, changing the adjustment coefficient is equivalent to rotating the descent direction, and by rotating the descent direction to an appropriate angle, a more optimal descent direction can be obtained.

Next, we will explain that when the adjustment coefficient takes on the special values of 0, 1, and $\frac{1}{2}$, it corresponds to algorithms with different orders of discretization error in the differential equations.

Consider the first order ordinary differential equation (gradient flow):
\begin{eqnarray}\label{eqn1-1}
	\left\{\begin{array}{l}
		\begin{aligned}
			\dot{x}(t)&=-\nabla f\left(x(t)\right), t\geq0,\\
			x(0)&=x^0.
		\end{aligned}
	\end{array}\right.
\end{eqnarray}
The equilibrium points $x^*=\lim_{t\rightarrow\infty} x(t)$ of (\ref{eqn1-1}) are the stationary points of unconstrained optimization problem (\ref{eqn1-2}), where $x(t)$ is the solution trajectory of (\ref{eqn1-1})\cite{AC20,CL21,SB16}. In general, the solution trajectory $x(t)$ is difficult to compute. Thus, some discrete methods are proposed to approximate $x(t)$.
The explicit Euler discrete method for problem (\ref{eqn1-1}) is related to the gradient descent algorithm \cite{Be17,LH20}
\begin{eqnarray*}
	x^{k+1}=x^k-h\nabla f(x^{k})
\end{eqnarray*}
for problem (\ref{eqn1-2}), and this corresponds to $\beta=0$ in Algorithm (\ref{eqn1-13}). The implicit Euler discrete method for problem (\ref{eqn1-1}) corresponds to the proximal point algorithm \cite{PB13,Ca19,CG22,JC21}
\begin{eqnarray*}\label{eqn1-6}
	x^{k+1}=x^k-h\nabla f(x^{k+1})
\end{eqnarray*}
for problem (\ref{eqn1-2}). Both the explicit Euler discrete method and the implicit Euler discrete method bring $O(h)$ global discrete error in terms of $\|x^k-x(t_k)\|$ \cite{Lin05}. Generally speaking, the proximal point algorithm is more stable than the gradient descent algorithm, but the proximal point algorithm involves implicit updates. If take the point obtained by the gradient descent algorithm as the approximation of the point $x^{k+1}$, then the proximal point algorithm can be approximated by extra-gradient algorithm (\ref{eqn1-71}). Notice that extra-gradient algorithm corresponds to $\beta=1$ in Algorithm (\ref{eqn1-13}).

The midpoint rectangle formula for (\ref{eqn1-1}) has higher global discrete error order of $O(h^2)$ \cite{Lin05}, and the discrete scheme is
\begin{eqnarray}\label{eqn1-8}
	x^{k+1}=x^k-h\nabla f\left(\frac{x^{k+1}+x^{k}}{2}\right).
\end{eqnarray}
Setting $z^k=\frac{x^{k+1}+x^{k}}{2}$, the iteration (\ref{eqn1-8}) can be rewritten as
\begin{eqnarray*}\label{eqn1-9}
	z^{k}&=&x^k-\frac{h}{2}\nabla f(z^{k}),\nonumber\\
	x^{k+1}&=&z^k+z^k-x^k,
\end{eqnarray*}
which corresponds to the extended proximal point algorithm for (\ref{eqn1-2}) \cite{PB13,HF09,Ne05}. As we all know, extended proximal point algorithm usually brings better performance than proximal point algorithm. It is important to note that in differential equations, the step size $h$ approaches zero to ensure that the iteration sequence $\{x^{k}\}$ approximate the solution trajectory $x(t)$, leading the limit of the iteration sequence to be the equilibrium point of the differential equation. However, in optimization problems, the primary concern is the final equilibrium point. To achieve faster convergence, a step size $h$ that does not approach zero is often used. This prevents the extrapolation coefficient $\alpha\in(0,1)$ in the extended proximal point algorithm: $x^{k+1}=z^k+\alpha (z^k-x^k)$ from reaching 1.

The trapezoidal formula for (\ref{eqn1-1}) also has global discrete error order of $O(h^2)$ \cite{Lin05}, and its discrete format is
\begin{eqnarray}\label{eqn1-10}
	x^{k+1}&=&x^k-\frac{h}{2}\left(\nabla f(x^{k+1})+\nabla f(x^{k})\right)\nonumber\\
	&=&x^k-h\left(\nabla f(x^{k})-\frac{1}{2}\left(\nabla f(x^{k})-\nabla f(x^{k+1})\right)\right).
\end{eqnarray}
However, iteration (\ref{eqn1-10}) is an implicit update. If we also take the point obtained by the gradient descent algorithm as the approximation of the point $x^{k+1}$, i.e.
\begin{eqnarray}\label{eqn1-11}
	z^k&=&x^k-h\nabla f(x^{k}),\nonumber\\
	x^{k+1}&=&x^k-h\left(\nabla f(x^{k})-\frac{1}{2}\left(\nabla f(x^{k})-\nabla f(z^{k})\right)\right).
\end{eqnarray}
(\ref{eqn1-11}) corresponds to Algorithm (\ref{eqn1-13}) with $\beta =\frac{1}{2}$. Since the trapezoidal formula has higher order of global discretization error than explicit and implicit Euler discretization, similar to the extended proximal point algorithm, Algorithm (\ref{eqn1-13}) may have better numerical performance when $\beta$ approaches to $\frac{1}{2}$, compared with when $\beta$ approaches to 0 or 1.

\section{Improved prediction correction algorithms with adjustment coefficient}

This section include three subsections. In subsection 3.1, we present two improved prediction-correction algorithms with adjustment coefficients for pseudo-convex objective functions under both constant and adaptive step sizes. In subsection 3.2, for convex objective functions, we introduce another improved prediction-correction algorithm with a wider range of adjustment coefficient. Additionally, the correction step can have a different step size from the prediction step. In subsection 3.3, we explain that when the objective function is convex and the adjustment coefficient is $\frac{1}{2}$, it can fully corresponds to the algorithm in (\ref{eqn1-11}).

\subsection{Improved prediction correction algorithms: pseudo-convex case}

In this subsection, we propose the improved prediction correction algorithms that use constant step size (IPC-C) and adaptive step size (IPC-A) for pseudo-convex optimization problem (\ref{eqn1-2}). To establish the convergence of the sequence $\{x^k\}$ generated by Algorithm IPC-C and IPC-A, we make some required assumptions firstly in the following.
\begin{assumption}\label{as1}
	\begin{itemize}
		\item[(i)] The stationary point set $X^{*}$ of problem $(\ref{eqn1-2})$ is non-empty.
		
		\item[(ii)] $f$ is pseudo-convex and $\nabla f$ is L-Lipschitz continuous.
		
	\end{itemize}
\end{assumption}

Algorithm 1 considers the iteration format (\ref{eqn1-13}) and uses a constant step size whose range related to the Lipschitz constant of $\nabla f$.

\begin{framed}
	\noindent\textbf{Algorithm 1 (pseudo-convex-IPC-C):}
	\begin{enumerate}[\bf Step 1:]
		\setcounter{enumi}{-1}
		\item Initialize $x^0\in \mathbb{R}^{n}$, $\epsilon>0$, $k=0$.
		\item If $\|\nabla f(x^{k})\|<\epsilon$, stop. Else go to step 2.
		\item Generate new iteration point.
		\begin{eqnarray}\label{eqn3-2-1}
			z^k&=&x^k-h\nabla f(x^{k}),\\\nonumber
			x^{k+1}&=&x^k-h\left(\nabla f(x^{k})-\beta(\nabla f(x^{k})-\nabla f(z^{k}))\right),
		\end{eqnarray}
		where $0<h<\frac{1}{L}$, $\beta\in\left(\frac{1-\sqrt{1-h^2L^2}}{h^2L^2},1\right]$.
		
		$k=k+1$. Go to Step 1.
	\end{enumerate}
\end{framed}

The Fej\'{e}r monotonicity of the sequence $\{x^k\}$ generated by Algorithm 1 with respect to the set $X^*$ is given below.

\begin{proposition}\label{prop3-2-1}
	Suppose that Assumption \ref{as1} holds, then the sequence $\{x^k\}$ generated by Algorithm 1 satisfies
	\begin{eqnarray*}\label{eqn3-2-3}
		\|x^{k+1}-x^*\|^2\leq\|x^{k}-x^*\|^2-\kappa_1\|x^{k}-z^k\|^2,
	\end{eqnarray*}
	where $\kappa_1=2\beta-1-\beta^2h^2L^2>0$.
\end{proposition}
\noindent{\bf\em Proof}\;~By the identity transformation, we obtain
\begin{eqnarray*}
	\|x^{k+1}-x^{*}\|^{2}
	&=&\|x^{k}-x^{*}\|^{2}+\|x^{k+1}-x^{k}\|^2+2\langle x^{k+1}-x^{k},x^{k}-x^{*}\rangle\nonumber\\
	&=&\|x^{k}-x^{*}\|^{2}+\|x^{k+1}-x^{k}\|^{2}+2\langle -h\nabla f(x^{k})+h\beta\left(\nabla f(x^k)-\nabla f(z^{k})\right),x^{k}-x^{*}\rangle\\
	&=&\|x^{k}-x^{*}\|^{2}+\|x^{k+1}-x^{k}\|^{2}-2\langle h\beta\nabla f(z^{k}),z^{k}-x^{*}\rangle\\
	&&-2\langle h\beta\nabla f(z^{k}),x^{k}-z^{k}\rangle+2\langle h(\beta-1)\nabla f(x^{k}),x^{k}-x^{*}\rangle.\nonumber
\end{eqnarray*}
Since $\nabla f$ is pseudo-monotone, $\nabla f(x^*)=0$ and $\beta\leq1$, we have
\begin{eqnarray}\label{eqn3-2-7}
	\|x^{k+1}-x^{*}\|^{2}
	&\leq&\|x^{k}-x^{*}\|^{2}+\|x^{k+1}-x^{k}\|^{2}-2\langle h\beta\nabla f(z^{k}),x^{k}-z^{k}\rangle.
\end{eqnarray}
Combining $(\ref{eqn3-2-1})$ and (\ref{eqn3-2-7}), it gets that
\begin{eqnarray}\label{eqn3-2-4}
	\|x^{k+1}-x^{*}\|^{2}
	&\leq&\|x^{k}-x^{*}\|^{2}+\|z^k-x^k+h\beta\nabla f(x^{k})-h\beta\nabla f(z^{k})\|^{2}\nonumber\\
	&&+2\beta\langle z^k-x^k+h\nabla f(x^{k})-h\nabla f(z^{k}),x^{k}-z^{k}\rangle\nonumber\\
	&=&\|x^{k}-x^{*}\|^{2}+\|z^k-x^k\|^2+\beta^2h^2\|\nabla f(x^{k})-\nabla f(z^{k})\|^2\nonumber\\
	&&+2\beta\langle h\nabla f(x^{k})-h\nabla f(z^{k}),z^k-x^k\rangle+2\beta\langle z^k-x^k+h\nabla f(x^{k})-h\nabla f(z^{k}),x^{k}-z^{k}\rangle\nonumber\\
	&=&\|x^{k}-x^{*}\|^{2}+\|z^k-x^k\|^2+\beta^2h^2\|\nabla f(x^{k})-\nabla f(z^{k})\|^2-2\beta\|x^{k}-z^{k}\|^2.
\end{eqnarray}
According to the $L$-Lipschitz continuity of $\nabla f$, (\ref{eqn3-2-4}) can be written as
\begin{eqnarray*}
	\|x^{k+1}-x^{*}\|^{2}
	&\leq&\|x^{k}-x^{*}\|^{2}+\|z^k-x^k\|^2+\beta^2h^2L^2\|x^{k}-z^{k}\|^2-2\beta\|x^{k}-z^{k}\|^2\nonumber\\
	&=&\|x^{k}-x^{*}\|^{2}-(2\beta-1-\beta^2h^2L^2)\|z^k-x^k\|^2.
\end{eqnarray*}
Due to $0<h<\frac{1}{L}$, $\beta\in\left(\frac{1-\sqrt{1-h^2L^2}}{h^2L^2},1\right]$, we have $2\beta-1-\beta^2h^2L^2>0$. Thus, the desired conclusion is proved. \qed

\begin{remark}\label{remark1}
	When $h\rightarrow0$, the lower bound of $\beta=\frac{1-\sqrt{1-h^2L^2}}{h^2L^2}\rightarrow\frac{1}{2}$. When $\beta=\frac{1}{2}$, Algorithm 1 becomes Algorithm (\ref{eqn1-11}). In the following numerical experiment, we find that, for a given $h$, the closer $\beta$ is to $\frac{1}{2}$, the better the numerical performance.
\end{remark}

The convergence of Algorithm 1 is given by the following theorem.
\begin{theorem}\label{the3-1}
	Suppose that Assumption \ref{as1} holds, then the sequence $\{x^k\}$ generated by Algorithm 1 converges to a stationary point of problem $(\ref{eqn1-2})$.
\end{theorem}
\noindent{\bf\em Proof}\;~Under Assumption \ref{as1}, the conclusion of Proposition \ref{prop3-2-1} holds. Thus
\begin{eqnarray*}
	h^2\kappa_1\sum_{k=0}^{\infty}\|\nabla f(x^k)\|^2=\kappa_1\sum_{k=0}^{\infty}\|x^{k}-z^k\|^2<\infty,
\end{eqnarray*}
and $\{x^k\}$ is bounded. Consequently, $\|\nabla f(x^{k})\|\rightarrow 0$, and there exists a subsequence $\{x^{k_j}\}$ converges to a cluster point $x^\infty$. By the continuity of the operator $\nabla f$ and $\|\cdot\|$, we have $\|\nabla f(x^{k_j})\|\rightarrow \|\nabla f(x^\infty)\|$. Combining with $\|\nabla f(x^k)\|\rightarrow 0$, we get $\|\nabla f(x^\infty)\|=0$. Thus $x^\infty\in X^*$. We next show that $\{x^k\}$ has a unique cluster point.

Since the sequence $\{\|x^k-x^*\|^2\}$ is monotonically decreasing and bounded below, thus $\{\|x^k-x^*\|^2\}$ converges. Taking $x^*=x^\infty$, we obtain $\{\|x^k-x^\infty\|^2\}$ has a limit. Taking the subsequence $\{x^{k_j}\}$, we have $\{\|x^{k_j}-x^\infty\|\}\rightarrow 0$. Thus the whole sequence $\{\|x^{k}-x^\infty\|\}\rightarrow 0$, i.e., $\{x^{k}\}\rightarrow x^\infty\in X^*$. \qed

The following theorem provides the convergence rate of $\|\nabla f(x^k)\|^2$ in the ergodic sense.

\begin{theorem}\label{the3-2}
	Suppose that Assumption \ref{as1} holds, then the sequence $\{x^k\}$ generated by Algorithm 1 satisfies
	\begin{eqnarray*}
		\frac{1}{K}\sum_{k=0}^{K-1}\|\nabla f(x^k)\|^2\leq\frac{\|x^0-x^*\|^2}{K\kappa_1h^2}.
	\end{eqnarray*}
\end{theorem}
\noindent{\bf\em Proof}\;~Under Assumption \ref{as1}, the conclusion of Proposition \ref{prop3-2-1} holds. Thus
	\begin{eqnarray*}
	\kappa_1 h^2\|\nabla f(x^k)\|^2=\kappa_1\|x^{k}-z^k\|^2\leq\|x^{k}-x^*\|^2-\|x^{k+1}-x^*\|^2.
\end{eqnarray*}
Sum the above inequality from 0 to $K-1$ to get
	\begin{eqnarray*}
	\kappa_1 h^2\sum_{k=0}^{K-1}\|\nabla f(x^k)\|^2\leq\|x^{0}-x^*\|^2.
\end{eqnarray*}
This means
	\begin{eqnarray*}
	\frac{1}{K}\sum_{k=0}^{K-1}\|\nabla f(x^k)\|^2\leq\frac{\|x^{0}-x^*\|^2}{K\kappa_1 h^2}.
\end{eqnarray*}
This proves the desired conclusion. \qed

The range of the step size $h$ and the adjustment coefficient $\beta$ in Algorithm 1 are determined by the global Lipschitz constant $L$ of $\nabla f$. However, generally speaking, $L$ is not easy to obtain or is difficult to estimate accurately. Even if $L$ is estimated well, it is only the global Lipschitz constant, which is much larger than the local Lipschitz constant. This will result in a small step size and slow convergence. Therefore, we present the following Algorithm IPC-A that does not rely on the global Lipschitz constant.
\begin{framed}
	\noindent\textbf{Algorithm 2 (pseudo-convex-IPC-A):}
	\begin{enumerate}[\bf Step 1:]
		\setcounter{enumi}{-1}
		\item  Initialize $x^0\in \mathbb{R}^{n}$, $0<\mu<\nu<1$, $\underline{ h}<1\leq\gamma_0^0\leq\overline{ h}$, $\theta\in(0,1)$, $\tau>1$, $\epsilon>0$, $k=0$.
		\item  If $\|\nabla f(x^{k})\|<\epsilon$, stop. Else go to step 2.
		\item Line search for $ h_k$.
		
		$l=0$;
		
		If $r_k(\gamma_l^k):=\gamma_l^k\|\nabla f(z^k(\gamma_l^k))-\nabla f(x^{k})\|/\|z^k(\gamma_l^k)-x^k\|>\nu$,
		
		\quad where $z^k(\gamma_l^k)=x^k-\gamma_l^k\nabla f(x^{k})$.
		
		\quad $\gamma_{l+1}^k=\gamma_l^k\theta*\min\{1,1/r_k(\gamma_l^k)\}$, $l=l+1$.
		
		Else
		
		\quad $ h_k=\gamma_l^k$.
		\item Generate new iteration point.
		\begin{eqnarray}\label{eqn3-2}
			z^k&=&x^k- h_k\nabla f(x^{k}),\\\nonumber
			x^{k+1}&=&x^k-h_k\left(\nabla f(x^{k})-\beta(\nabla f(x^{k})-\nabla f(z^{k}))\right),
		\end{eqnarray}
		\quad where $\beta\in\left(\frac{1-\sqrt{1-\nu^2}}{\nu^2},1\right]$.
		\item Adaptively select $\gamma_0^{k+1}$, the initial line search step size for iteration $k+1$.
		
		If $r_k( h_k)\leq\mu$, then $\gamma_0^{k+1}=\tau h_k$.
		
		Else $\gamma_0^{k+1}= h_k$.
		
		$\gamma_0^{k+1}=\textbf{P}_{[\underline{ h},\overline{ h}]}(\gamma_0^{k+1})$.
		
		$k=k+1$. Go to Step 1.
	\end{enumerate}
\end{framed}

We first illustrate by two lemmas that the line search in Algorithm 2 finitely terminates and the step size $h_k$ has a consistent lower bound.

\begin{lemma}\label{lem3-1-1}
	Suppose that $\nabla f$ is L-Lipschitz continuous, then the line search for $ h_k$ in Algorithm 2 terminates in finite steps.
\end{lemma}
\noindent{\bf\em Proof}\;~If there exists $k$ such that the line search can not finitely terminated, then for all $l>0$,
\begin{eqnarray*}\label{eqn3-1-1}
	\gamma_l^k\|\nabla f(z^k(\gamma_l^k))-\nabla f(x^{k})\|/\|z^k(\gamma_l^k)-x^k\|>\nu.
\end{eqnarray*}
By the $L$-Lipschitz continuity of $\nabla f$, we have
\begin{eqnarray*}\label{eqn3-1-2}
	\gamma_l^kL\geq\gamma_l^k\|\nabla f(z^k(\gamma_l^k))-\nabla f(x^{k})\|/\|z^k(\gamma_l^k)-x^k\|>\nu.
\end{eqnarray*}
Thus, $\gamma_l^k>\frac{\nu}{L}$, $\forall l>0$. Letting $l\rightarrow\infty$, we obtain $\frac{\nu}{L}\leq0$. This contradicts with $0<\nu<1$. \qed

\begin{lemma}\label{lem3-1-2}
	Suppose that $\nabla f$ is L-Lipschitz continuous, then the step size $ h_k$ in Algorithm 2 satisfies $ h_k\geq h_{\min}$, where $ h_{\min}=\min\left\{\underline{ h},\frac{\nu\theta}{\overline{ h}L^2}\right\}$.
\end{lemma}
\noindent{\bf\em Proof}\;~If $ h_k=\gamma_0^k$, then $ h_k\geq\underline{ h}$. If $ h_k=\gamma_l^k, l\geq 1$, then
\begin{eqnarray*}
	\gamma_{l-1}^k\|\nabla f(z^k(\gamma_{l-1}^k))-\nabla f(x^{k})\|/\|z^k(\gamma_{l-1}^k)-x^k\|>\nu.
\end{eqnarray*}
By the $L$-Lipschitz continuity of $\nabla f$, we get
\begin{eqnarray*}
	\gamma_{l-1}^kL\geq\gamma_{l-1}^k\|\nabla f(z^k(\gamma_{l-1}^k))-\nabla f(x^{k})\|/\|z^k(\gamma_{l-1}^k)-x^k\|>\nu.
\end{eqnarray*}
That is $\gamma_{l-1}^k>\frac{\nu}{L}$. Multiplying by $\theta*\min\left\{1,1/r_k(\gamma_{l-1}^k)\right\}$ at both sides of this inequality, we have
\begin{eqnarray}\label{eqn3-1-3}
	h_k=\gamma_l^k>\frac{\nu\theta}{L}*\min\left\{1,1/r_k(\gamma_{l-1}^k)\right\}.
\end{eqnarray}
Due to
\begin{eqnarray*}
	r_k(\gamma_{l-1}^k)&=&\gamma_{l-1}^k\|\nabla f(z^k(\gamma_{l-1}^k))-\nabla f(x^{k})\|/\|z^k(\gamma_{l-1}^k)-x^k\|\\
	&\leq&\frac{\gamma_{l-1}^kL\|z^k(\gamma_{l-1}^k)-x^k\|}{\|z^k(\gamma_{l-1}^k)-x^k\|}=\gamma_{l-1}^kL\leq\gamma_{0}^kL\leq\overline{ h }L,
\end{eqnarray*}
we obtain
\begin{eqnarray}\label{eqn3-1-4}
	\frac{1}{r_k(\gamma_{l-1}^k)}\geq\frac{1}{\overline{ h}L}.
\end{eqnarray}
From (\ref{eqn3-1-3}) and (\ref{eqn3-1-4}), we have $ h_k>\frac{\nu\theta}{L}$ or $ h_k>\frac{\nu\theta}{L}\cdot\frac{1}{\overline{ h}L}$. We can always make $L\geq 1$. Combining with $\overline{ h}\geq1$, we have $\frac{1}{\overline{ h}L}\leq1$. Thus, $ h_k>\frac{\nu\theta}{L}\cdot\frac{1}{\overline{ h}L}=\frac{\nu\theta}{\overline{ h}L^2}$.

Combining the above two situations, we have $ h_k\geq h_{\min}=\min\left\{\underline{h},\nu\theta/\overline{ h}L^2\right\}$. \qed

Similar to Proposition \ref{prop3-2-1}, The sequence $\{x^k\}$ generated by Algorithm 2 is Fej\'{e}r monotonicity with respect to the set $X^*$.

\begin{proposition}\label{prop3-1}
	Suppose that Assumption \ref{as1} holds, then the sequence $\{x^k\}$ generated by Algorithm 2 satisfy
	\begin{eqnarray*}\label{eqn3-3}
		\|x^{k+1}-x^*\|^2\leq\|x^{k}-x^*\|^2-\kappa_2\|x^{k}-z^k\|^2,
	\end{eqnarray*}
	where $\kappa_2=2\beta-1-\beta^2\nu^2>0$.
\end{proposition}
\noindent{\bf\em Proof}\;~Consistent with the first half of the proof in Proposition \ref{prop3-2-1},
the inequality (\ref{eqn3-2-4}) is obtained. Combining (\ref{eqn3-2-4}) with the choice of $ h_k$ in Algorithm 2, we get
\begin{eqnarray*}\label{eqn3-11}
	\|x^{k+1}-x^{*}\|^{2}
	&\leq&\|x^{k}-x^{*}\|^{2}+\|z^k-x^k\|^2+\beta^2\nu^2\|x^{k}-z^{k}\|^2-2\beta\|x^{k}-z^{k}\|^2\nonumber\\
	&=&\|x^{k}-x^{*}\|^{2}-(2\beta-1-\beta^2\nu^2)\|z^k-x^k\|^2.
\end{eqnarray*}
Due to $0<\nu<1$, $\beta\in\left(\frac{1-\sqrt{1-\nu^2}}{\nu^2},1\right]$, we have $2\beta-1-\beta^2\nu^2>0$. Thus, the desired conclusion is obtained. \qed

The convergence of Algorithm 2 is given by the following theorem.
\begin{theorem}\label{the3-3}
	Suppose that Assumption \ref{as1} holds, then the sequence $\{x^k\}$ generated by Algorithm 2 converges to a stationary point of problem $(\ref{eqn1-2})$.
\end{theorem}
\noindent{\bf\em Proof}\;~Under Assumption \ref{as1}, the conclusion of Proposition \ref{prop3-1} holds. Thus
\begin{eqnarray*}
	h_{\min}^2\kappa_2\sum_{k=0}^{\infty}\|\nabla f(x^k)\|^2\leq h_k^2\kappa_2\sum_{k=0}^{\infty}\|\nabla f(x^k)\|^2=\kappa_2\sum_{k=0}^{\infty}\|x^{k}-z^k\|^2<\infty,
\end{eqnarray*}
and $\{x^k\}$ is bounded. Consequently, $\|\nabla f(x^{k})\|\rightarrow 0$, and there exists a subsequence $\{x^{k_j}\}$ converges to a cluster point $x^\infty$. The latter proof is consistent with Theorem \ref{the3-1}. \qed

Consistent with the proof of Theorem \ref{the3-2} and combined with Lemma \ref{lem3-1-2}, we can conclude the following convergence rate result.

\begin{theorem}\label{the3-5}
	Suppose that Assumption \ref{as1} holds, then the sequence $\{x^k\}$ generated by Algorithm 2 satisfies
	\begin{eqnarray*}
		\frac{1}{K}\sum_{k=0}^{K-1}\|\nabla f(x^k)\|^2\leq\frac{\|x^0-x^*\|^2}{K\kappa_2h_{\min}^2}.
	\end{eqnarray*}
\end{theorem}

\subsection{Improved prediction correction algorithm: convex case}

In this subsection, we present another improved prediction-correction algorithm for the convex optimization problem (\ref{eqn1-2}), in which a wider range of adjustment coefficient can be used, and the correction step uses a different step size from the prediction step. Then, we show that if the step size in prediction step takes a special value, the optimal step size in the correction step of has a lower bound of $\frac{1}{2}$. Consequently, the algorithm can corresponds exactly to Algorithm (\ref{eqn1-11}).

\begin{framed}
	\noindent\textbf{Algorithm 3 (convex-IPC):}
	\begin{enumerate}[\bf Step 1:]
		\setcounter{enumi}{-1}
		\item  Initialize $x^0\in \mathbb{R}^{n}$, $0<\mu<\nu<1$, $0<\underline{ h}<1\leq\gamma_0^0\leq\overline{ h}<\frac{4}{L}$, $0\leq\beta \leq1$, $\theta\in(0,1)$, $\tau>1$, $\eta\in(0,2),$ $\epsilon>0$, $k=0$.
		\item  If $\|\nabla f(x^k)\|<\epsilon$, stop. Else go to step 2.
		\item Line search for $ h_k$.
		
		$l=0$;
		
		If $r_k(\gamma_l^k):=\gamma_l^k\|\nabla f(z^k(\gamma_l^k))-\nabla f(x^{k})\|/\|z^k(\gamma_l^k)-x^k\|>\nu$,
		
		\quad where $z^k(\gamma_l^k)=x^k-\gamma_l^k\nabla f(x^{k})$.
		
		\quad $\gamma_{l+1}^k=\gamma_l^k\theta*\min\{1,1/r_k(\gamma_l^k)\}$, $l=l+1$.
		
		Else
		
		\quad $ h_k=\gamma_l^k$.
		\item Generate new iteration point.
		\begin{eqnarray}
			z^k&=&x^k- h_k\nabla f(x^{k}),\label{eqn3-16}\\
			x^{k+1}&=&x^k-\eta\alpha_{k}h_k\left(\nabla f(x^k)-\beta(\nabla f(x^k)-\nabla f(z^k))\right),\label{eqn3-17}
		\end{eqnarray}
		where $\alpha_{k}=\frac{(1-\beta )\left(1-\frac{Lh_k}{4}\right)\|x^k-z^k\|^2+\beta \langle x^k-z^k,h_k\nabla f(z^k)\rangle}{h_k^2\|\nabla f(x^k)-\beta(\nabla f(x^k)-\nabla f(z^k))\|^2}$.
		\item Adaptively select $\gamma_0^{k+1}$, the initial line search step size for iteration $k+1$.
		
		If $r_k( h_k)\leq\mu$, then $\gamma_0^{k+1}=\tau h_k$.
		
		Else $\gamma_0^{k+1}= h_k$.
		
		$\gamma_0^{k+1}=\textbf{P}_{[\underline{ h},\overline{ h}]}(\gamma_0^{k+1})$.
		
		$k=k+1$. Go to Step 1.
	\end{enumerate}
\end{framed}

\begin{remark}\label{remark2}
In Algorithm 3, the calculation of $\alpha_k$ requires the global Lipschitz constant $L$, but we still use adaptive line search in Step 2. This is because the line search seeks a smaller local Lipschitz constant, allowing for a larger step size $h_k$, and the algorithm usually converges faster.
\end{remark}

We first give assumptions required by Algorithm 3.

\begin{assumption}\label{as2}
	\begin{itemize}
		\item[(i)] The solution set $X^*$ of problem $(\ref{eqn1-2})$ is non-empty.
		
		\item[(ii)] $f$ is convex and  $\nabla f$ is $L$-Lipschitz continuous.
	\end{itemize}
\end{assumption}

We next demonstrate below through a lemma that the direction in (\ref{eqn3-17}) is a descent direction for the distance function $\frac{1}{2}\|x^k-x^*\|^2$.

\begin{lemma}\label{lem3-4}
	Suppose Assumption \ref{as2} holds. Let $\{x^k\}, \{z^k\}$ are generated by Algorithm 3, and define  $$d_k:=h_k\left(\nabla f(x^k)-\beta (\nabla f(x^k)-\nabla f(z^k))\right),$$
	where $0\leq\beta \leq1$. Then
	\begin{eqnarray}
		\langle x^k-x^*,d^k\rangle&\geq& (1-\beta )\left(1-\frac{Lh_k}{4}\right)\|x^k-z^k\|^2+\beta \langle x^k-z^k,h_k\nabla f(z^k)\rangle\label{eqn3-14}\\
		&\geq&\kappa_3\|x^k-z^k\|^2\label{eqn3-36},
	\end{eqnarray}
	where $\kappa_3=(1-\beta )\left(1-\frac{Lh_k}{4}\right)+\beta (1-\nu)>0.$
\end{lemma}
\noindent{\bf\em Proof}\; According to the definition of $d^k$, we have
\begin{eqnarray}\label{eqn3-35}
	\langle x^k-x^*,d^k\rangle&=& \langle x^k-x^*,h_k\nabla f(x^k)-\beta h_k(\nabla f(x^k)-\nabla f(z^k))\rangle.
\end{eqnarray}
On the one hand, due to
\begin{eqnarray*}\label{eqn3-6}
	x^k-z^k=h_k\nabla f(x^k)=h_k\nabla f(x^k)-h_k\nabla f(x^*),
\end{eqnarray*}
we can obtain
\begin{eqnarray*}\label{eqn3-7}
	\langle x^k-x^*, x^k-z^k\rangle&=&\|x^k-z^k\|^2+h_k\langle z^k-x^*,\nabla f(x^k)-\nabla f(x^*)\rangle\nonumber\\
	&=&\|x^k-z^k\|^2+h_k\langle x^k-x^*,\nabla f(x^k)-\nabla f(x^*)\rangle-h_k\langle x^k-z^k,\nabla f(x^k)-\nabla f(x^*)\rangle.
\end{eqnarray*}
According to Lemma \ref{lem2-3} and the inequality $2ab\leq a^2+b^2$ with $a=\frac{\sqrt{Lh_k}}{2}(x^k-z^k)$, $b=\sqrt{\frac{h_k}{L}}(\nabla f(x^k)-\nabla f(x^*))$, we have
\begin{eqnarray*}\label{eqn3-39}
	\langle x^k-x^*, x^k-z^k\rangle
	&\geq&\|x^k-z^k\|^2+h_k\frac{1}{L}\|\nabla f(x^k)-\nabla f(x^*)\|^2-h_k\frac{1}{L}\|\nabla f(x^k)-\nabla f(x^*)\|^2-\frac{Lh_k}{4}\|x^k-z^k\|^2\nonumber\\
	&=&\left(1-\frac{Lh_k}{4}\right)\|x^k-z^k\|^2.
\end{eqnarray*}
This means
\begin{eqnarray}\label{eqn3-8}
	\langle x^k-x^*, h_k\nabla f(x^k)\rangle
	&\geq&\left(1-\frac{Lh_k}{4}\right)\|x^k-z^k\|^2.
\end{eqnarray}
On the other hand,
\begin{eqnarray}\label{eqn3-310}
	\langle x^k-x^*, h_k\nabla f(z^k)\rangle
	&=&\langle z^k-x^*,h_k\nabla f(z^k)\rangle +\langle x^k-z^k, h_k\nabla f(z^k)\rangle.
\end{eqnarray}
Combining this with $\langle z^k-x^*,h_k\nabla f(x^*)\rangle=0$ and the monotonicity of the operator $\nabla f$, we obtain
\begin{eqnarray}\label{eqn3-320}
	\langle z^k-x^*,h_k\nabla f(z^k)\rangle\geq0.
\end{eqnarray}
Combining (\ref{eqn3-35}), (\ref{eqn3-8}), (\ref{eqn3-310}), (\ref{eqn3-320}) and $0\leq\beta \leq1$, we obtain (\ref{eqn3-14}).

According to the iteration format (\ref{eqn3-16}), Cauchy-Schwarz inequality and the line search condition in Algorithm 3, we have
\begin{eqnarray}\label{eqn3-33}
	\langle x^k-z^k, h_k\nabla f(z^k)\rangle&=&\langle x^k-z^k, x^k-z^k-h_k\nabla f(x^k)+h_k\nabla f(z^k)\rangle\nonumber\\
	&\geq&\|x^k-z^k\|^2-h_k\|x^k-z^k\|\|\nabla f(x^k)-\nabla f(z^k)\|\nonumber\\
	&\geq&(1-\nu)\|x^k-z^k\|^2.
\end{eqnarray}
Combining (\ref{eqn3-310}), (\ref{eqn3-320}), (\ref{eqn3-33}), we obtain
\begin{eqnarray}\label{eqn3-34}
	\langle x^k-x^*, h_k\nabla f(z^k)\rangle&\geq&(1-\nu)\|x^k-z^k\|^2.
\end{eqnarray}
The inequality (\ref{eqn3-36}) can be obtained by combining the inequalities (\ref{eqn3-35}), (\ref{eqn3-8}), (\ref{eqn3-34}) and $0\leq\beta \leq1$.
\qed

We next show that $\alpha_k$ has a consistent lower bound.

\begin{lemma}\label{lem3-5}
	Under the assumption that the operator $\nabla f$ is $L$-Lipschitz continuous, the step size $\alpha_k$ in Algorithm 3 satisfies $\alpha_k\geq \alpha_{\min}$, where
	$$\alpha_{\min}=\frac{\left((1-\beta )\left(1-\frac{L\overline{h}}{4}\right)+\beta (1-\nu)\right)}{(2+2\beta ^2\nu^2)}.$$
\end{lemma}
\noindent{\bf\em Proof}\; According to (\ref{eqn3-33}) and the definition of $\alpha_k$, we have
\begin{eqnarray*}\label{eqn3-21}
	\alpha_{k}&=&\frac{(1-\beta )\left(1-\frac{Lh_k}{4}\right)\|x^k-z^k\|^2+\beta \langle x^k-z^k,h_k\nabla f(z^k)\rangle}{h_k^2\|\nabla f(x^k)-\beta(\nabla f(x^k)-\nabla f(z^k))\|^2}\\
	&\geq&\frac{\left((1-\beta )\left(1-\frac{Lh_k}{4}\right)+\beta (1-\nu)\right)\|x^k-z^k\|^2}{\|x^k-z^k-\beta h_k\left(\nabla f(x^k)-\nabla f(z^k)\right)\|^2}\\
	&\geq&\frac{\left((1-\beta )\left(1-\frac{Lh_k}{4}\right)+\beta (1-\nu)\right)\|x^k-z^k\|^2}{2\|x^k-z^k\|^2+2\|\beta h_k\left(\nabla f(x^k)-\nabla f(z^k)\right)\|^2}.
\end{eqnarray*}
Due to Step 2 in Algorithm 3, we have
\begin{eqnarray*}\label{eqn3-21}
	\alpha_{k}
	&\geq&\frac{\left((1-\beta )\left(1-\frac{Lh_k}{4}\right)+\beta (1-\nu)\right)\|x^k-z^k\|^2}{(2+2\beta ^2\nu^2)\|x^k-z^k\|^2}\\
	&=&\frac{\left((1-\beta )\left(1-\frac{Lh_k}{4}\right)+\beta (1-\nu)\right)}{(2+2\beta ^2\nu^2)}\\
	&\geq&\frac{\left((1-\beta )\left(1-\frac{L\overline{h}}{4}\right)+\beta (1-\nu)\right)}{(2+2\beta ^2\nu^2)}:=\alpha_{\min}.
\end{eqnarray*}
\qed

The following proposition states that the sequence $\{x^k\}$ generated by Algorithm 3 is Fej\'{e}r monotone with respect the solution set $X^*$.

\begin{proposition}\label{prop3-2}
	Under Assumption \ref{as2}, the sequence $\{x^k\}$ generated by Algorithm 3 satisfies
	\begin{eqnarray*}\label{eqn3-18}
		\|x^{k+1}-x^*\|^2\leq\|x^{k}-x^*\|^2-\kappa_4\|x^{k}-z^k\|^2,
	\end{eqnarray*}
	where $\kappa_4=\eta(2-\eta)\alpha_{\min}\left((1-\beta )\left(1-\frac{L\overline{h}}{4}\right)+\beta (1-\nu)\right)>0$.
\end{proposition}
\noindent{\bf\em Proof}\; From the iteration format (\ref{eqn3-17}), we have
\begin{eqnarray*}\label{eqn3-19}
	\|x^{k+1}-x^*\|^2=\|x^{k}-x^*\|^2-2\eta\alpha_{k}\langle x^k-x^*,d^k\rangle+\eta^2\alpha_{k}^2\|d^k\|^2.
\end{eqnarray*}
Combining (\ref{eqn3-14}), the choice of $\alpha_k$ in Step 3, and the above equality, we obtain
\begin{eqnarray*}\label{eqn3-20}
	&&\|x^{k+1}-x^*\|^2\\
	&\leq&\|x^{k}-x^*\|^2+\eta^2\alpha_{k}^2\|d^k\|^2-2\eta\alpha_{k}\left((1-\beta )\left(1-\frac{Lh_k}{4}\right)\|x^k-z^k\|^2+\beta \langle x^k-z^k,h_k\nabla f(z^k)\rangle\right)\\
	&=&\|x^{k}-x^*\|^2-\eta(2-\eta)\alpha_k\left((1-\beta )\left(1-\frac{Lh_k}{4}\right)\|x^k-z^k\|^2+\beta \langle x^k-z^k,h_k\nabla f(z^k)\rangle\right)\\
	&\leq&\|x^{k}-x^*\|^2-\eta(2-\eta)\alpha_{\min}\left((1-\beta )\left(1-\frac{L\overline{h}}{4}\right)+\beta (1-\nu)\right)\|x^k-z^k\|^2.\\
\end{eqnarray*}
where the equality is according to the setting of $\alpha_k$ in Step 3 in Algorithm 3, and the last inequality is due to (\ref{eqn3-36}), $\alpha_k\geq\alpha_{\min}$ and $h_k\leq\overline{h}$.
\qed

Similar to the proofs of Theorem \ref{the3-1} and Theorem \ref{the3-3}, we can derive the following convergence theorem.
\begin{theorem}\label{the3-4}
	Under Assumption \ref{as2}, the sequence $\{x^k\}$ generated by Algorithm 3 converges to a solution point of problem $(\ref{eqn1-2})$.
\end{theorem}

Similar to the proof of Theorem \ref{the3-5}, we can conclude the following convergence rate result.

\begin{theorem}\label{the3-6}
	Suppose that Assumption \ref{as2} holds, then the sequence $\{x^k\}$ generated by Algorithm 3 satisfies
	\begin{eqnarray*}
		\frac{1}{K}\sum_{k=0}^{K-1}\|\nabla f(x^k)\|^2\leq\frac{\|x^0-x^*\|^2}{K\kappa_4h_{\min}^2}.
	\end{eqnarray*}
\end{theorem}

\subsection{The relationship between Algorithm 3 and Algorithm (\ref{eqn1-11})}

If we take the step size $h_k\equiv \frac{0.9}{L}$, Algorithm 3 can be written as
\begin{eqnarray*}\label{eqn3-25}
	z^k&=&x^k-\frac{0.9}{L}\nabla f(x^{k}),\nonumber\\
	x^{k+1}&=&x^k-\eta\alpha_k \frac{0.9}{L}\left(\nabla f(x^{k})-\beta\left(\nabla f(x^{k})-\nabla f(z^{k})\right)\right).
\end{eqnarray*}
Based on the expression of $\alpha_k$ in Algorithm 3, after transformation, we obtain
\begin{eqnarray}\label{eqn3-22}
	\alpha_{k}&=&\frac{\frac{3.1}{4}(1-\beta )\|x^k-z^k\|^2}{\left\|x^k-z^k-\frac{0.9}{L}\beta\left(\nabla f(x^{k})-\nabla f(z^{k})\right)\right\|^2}+\frac{\beta \langle x^k-z^k,x^k-z^k-\frac{0.9}{L}\left(\nabla f(x^{k})-\nabla f(z^{k})\right)\rangle}{\left\|x^k-z^k-\frac{0.9}{L}\beta\left(\nabla f(x^{k})-\nabla f(z^{k})\right)\right\|^2}\nonumber\\
	&=&\frac{1}{2}\frac{\left(\frac{3.1+0.9\beta}{2}\right)\|x^k-z^k\|^2-\frac{1.8}{L}\beta \langle x^k-z^k,\nabla f(x^k)-\nabla f(z^k)\rangle}{\left\|\left(x^k-z^k-\frac{0.9}{L}\beta\left(\nabla f(x^{k})-\nabla f(z^{k})\right)\right)\right\|^2}\nonumber\\
	&=&\frac{1}{2}\left(1+\frac{\left(\frac{1.1+0.9\beta }{2}\right)\|x^k-z^k\|^2-\frac{0.81}{L^2}\beta ^2\|\nabla f(x^k)-\nabla f(z^k)\|^2}{\left\|x^k-z^k-\frac{0.9}{L}\beta\left(\nabla f(x^{k})-\nabla f(z^{k})\right)\right\|^2}\right).
\end{eqnarray}
Due to the Assumption that the operator $\nabla f$ is $L$-Lipschitz continuous, we have under $\beta \in[0,1],$
\begin{eqnarray}\label{eqn3-223}
	&&\left(\frac{1.1+0.9\beta }{2}\right)\|x^k-z^k\|^2-\frac{0.81}{L^2}\beta ^2\|\nabla f(x^k)-\nabla f(z^k)\|^2\nonumber\\
	&\geq&\left(\left(\frac{1.1+0.9\beta }{2}\right)-0.81\beta ^2\right)\|x^k-z^k\|^2>0.
\end{eqnarray}
Combining (\ref{eqn3-22}) and (\ref{eqn3-223}), we obtain $\alpha_k>\frac{1}{2}$. Thus, we can take $\eta=\frac{1}{\alpha_k}\in(0,2)$. So that $\eta\alpha_k=1$. At this point, Algorithm 3 can be written as
\begin{eqnarray*}
	z^k&=&x^k-\frac{0.9}{L}\nabla f(x^{k}),\nonumber\\
	x^{k+1}&=&x^k-\frac{0.9}{L}\left(\nabla f(x^{k})-\beta\left(\nabla f(x^{k})-\nabla f(z^{k})\right)\right).
\end{eqnarray*}
where $\beta\in[0,1]$, which is in the form of Algorithm (\ref{eqn1-11}). Notice that the range of $\beta$ in Algorithm 3 is wider than Algorithm 1 and 2.
If $\beta =0$, Algorithm 3 is the gradient descent algorithm. If $\beta =1$, Algorithm 3 is the extra-gradient algorithm. If $\beta =\frac{1}{2}$, Algorithm 3 becomes Algorithm (\ref{eqn1-11}) which is related to the trapezoidal formula.

\section{Numerical Experiment}

In this section, we provide some examples of both pseudo-convex and convex optimization problems to demonstrate the effectiveness of the algorithms proposed in this paper. The experiments are operated with Matlab R2021b on a Windows 11 with a 3.30 Ghz processor and 32 GB of memory.

\textbf{ \\Experiment 1: Pseudo-convex fractional programming problem}
\vspace{0.3cm}

Consider the fractional program problem of the form \cite{BM19}:
$$\min f(x)=\frac{G(x)}{h(x)},$$
where $G(x)$ is a convex quadratic function and $h(x)$ is a positive linear function. Apparently, such a problem is pseudo-convex and its gradient is pseudo-monotone \cite{Ma70}. Specifically, $G(x)=\dfrac{1}{2}x'Qx+c'x+q$, $h(x)=r'x+t$. $Q=MM'+I$, $M$ is a $n\times n$ matrix and each element is selected from a uniform distribution on the interval (0,1) randomly. $I$ is the $n\times n$ identity matrix. $r$ and $c$ are $n\times 1$ vectors drawn uniformly at random from $(0,2)$. $q$ is a random number in (1,2) and $t=1+4n$.

To avoid computing the Lipschitz constant of the gradient of the objective function in this experiment, we only provide the numerical results for Algorithm 2. The starting point $x^{0}$ is randomly chosen in $(1,10)^{n}$. The parameters in our Algorithm 2 are set as: $\mu=0.3, \nu=0.5, \tau=1.5, \theta=0.67, \gamma_0^0=1, \overline{h}=3, \underline{h}=10^{-6}, \epsilon=10^{-3}$. The stopping criterion is chosen as $\|\nabla f(x^{k})\|\leq\epsilon$.

\begin{figure}[htbp]
		\begin{center}
	\subfigure{
		\includegraphics[width=3.1in]{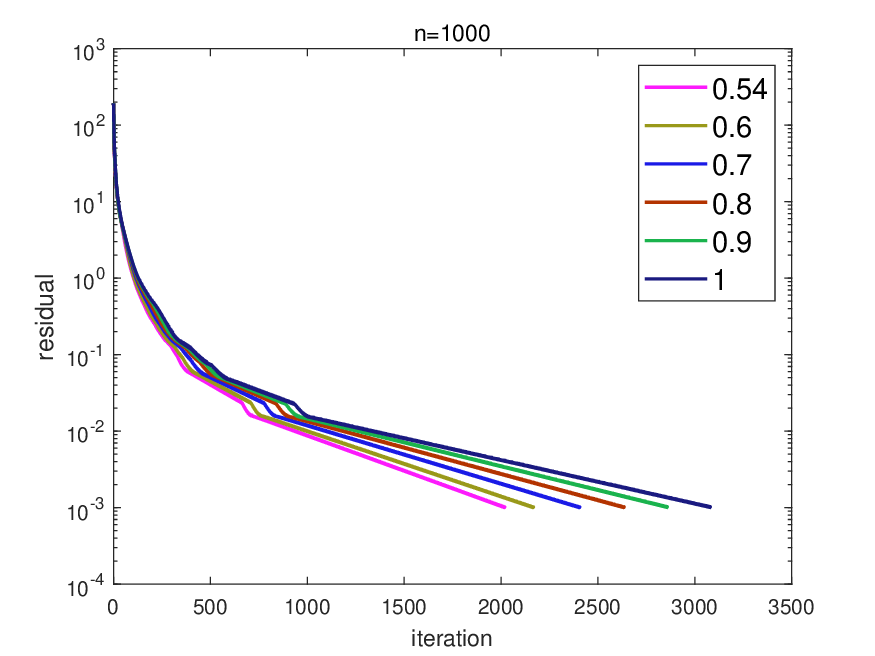}
	}
	\subfigure{
		\includegraphics[width=3.1in]{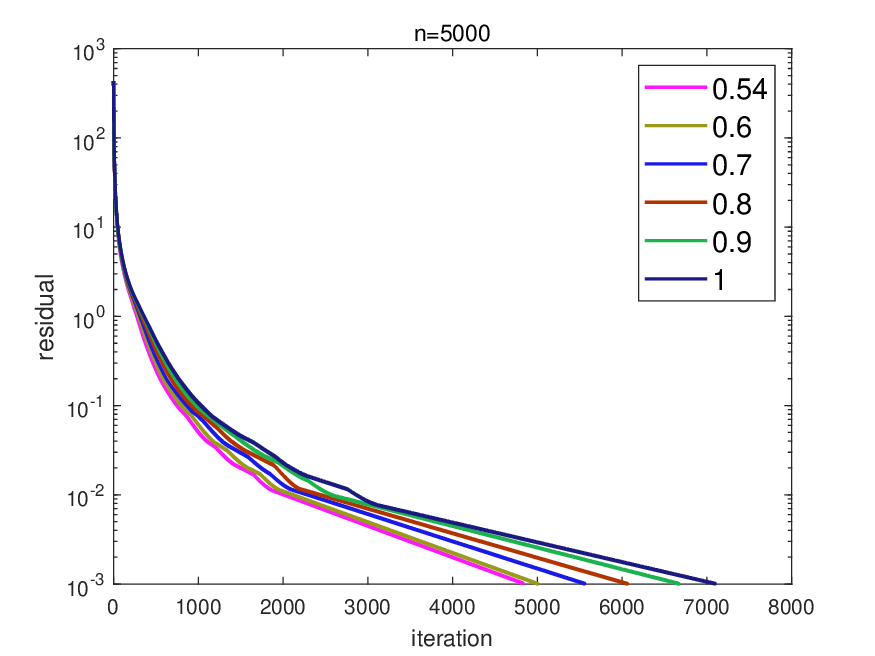}
	}
	\caption{Numerical performance of Algorithm 2 with different $\beta $ for Experiment 1}
		\label{fig1}
		\end{center}

\end{figure}

We carry out our experiments for $n\in\{1000,2000,5000\}$ under different $\beta$. Notice that $\beta=0.54$ is the minimum value that can be attained in Algorithm 2 when $\nu=0.5$. The number of iterations and CPU time under different dimensions and $\beta$ of Algorithm 2 are given in Table \ref{table1}. The residual $\log\left(\|\nabla f(x^{k})\|\right)$ vs. number of iterations is given in Figure \ref{fig1}.

\begin{center}
	\begin{table}[htbp]
		\centering
		\caption{Comparison of the efficiency of Algorithm 2 under different $\beta$ for Experiment 1}
		\begin{tabular}
			{m{1.2cm}<{\centering}|m{1.2cm}<{\centering}m{1.2cm}<{\centering}|m{1.2cm}<{\centering}m{1.2cm}<{\centering}|m{1.2cm}<{\centering}m{1.2cm}<{\centering}}
			\toprule
			&\multicolumn{2}{c|}{$n=1000$}&\multicolumn{2}{c|}{$n=2000$}&\multicolumn{2}{c}{$n=5000$}\\\hline
			&Iter.&CPU&Iter.&CPU&Iter.&CPU\\
			$\beta=0.54$&\textbf{2022}&\textbf{1.6364}&\textbf{3123}&\textbf{14.6207}&\textbf{4840}&\textbf{177.0001}\\
			$\beta=0.6$&2170&1.7270&3244&14.9632&5008&186.8468\\
			$\beta=0.7$&2409&2.0179&3616&16.6457&5561&201.4298\\
			$\beta=0.8$&2638&2.0183&3964&18.5336&6065&210.1741\\
			$\beta=0.9$&2863&2.2903&4383&19.9577&6675&223.3490\\
			$\beta=1$&3085&2.4040&4723&21.7784&7101&237.7850\\
			\bottomrule
		\end{tabular}
		\label{table1}
	\end{table}
\end{center}

\vspace {-1cm}

The comparison of numerical results of Algorithm 2 under different dimensions and $\beta$ are shown in Table \ref{table1} and Figure \ref{fig1}. When $\beta=1$, Algorithm 2 reduces to the extra-gradient algorithm. It can be seen from these tables and figures that at various dimensions, the numerical performance is best when $\beta=0.54$ in Algorithm 2, and improves by approximately 32\% compared to the extra-gradient algorithm. This optimal $\beta$ is close to 0.5, which is consistent with our previous explanation from the perspective of differential equations.

\textbf{ \\Experiment 2: Smooth convex optimization problem}
\vspace{0.3cm}

Consider the following differentiable convex optimization problem, in which we refer to the example in \cite{HL02}.
\begin{eqnarray*}\label{eqn4-29}
	&\min\limits_{x\in \mathbb{R}^n}&~f(x)=x^T\arctan x-\frac{1}{2}\sum_{i=1}^{n}\ln(1+x_i^2)+\frac{1}{2}x^T Mx+q^Tx.
\end{eqnarray*}
The matrix $M=A^TA+B$, where $A$ is an $n\times n$ matrix whose entries are randomly generated in the interval (-5,5) and a skew-symmetric matrix $B$ is generated in the same way. The vector $q$ is generated from a uniform distribution in the interval (-500,500).
The gradient of $f$ is $\nabla f(x)= \arctan x + Mx+q.$ The Lipschitz constant of $\nabla f$ is $L=\|M\|+1$.

We carry out our experiments under dimension $n\in \{1000, 2000, 5000\}$, and compare the numerical performance of Algorithm 3 under different $\beta$. Each algorithm comes from the same initial point $x^0$, whose entries are randomly generated in the interval (0,1). The parameters in our Algorithm 3 are set as: $\mu=0.4,$ $\nu=0.9,$ $\underline{h}=10^{-6}$, $\gamma_{0}^0=\overline{h}=2/L$, $\theta=0.7,$ $\tau=1.5$, $\eta=1.9$, $\epsilon=10^{-3}$. The stopping criterion is set as $\|\nabla f(x^k)\|<\epsilon$.

\begin{figure}[htbp]
			\begin{center}
	\subfigure{
		\includegraphics[width=3.1in]{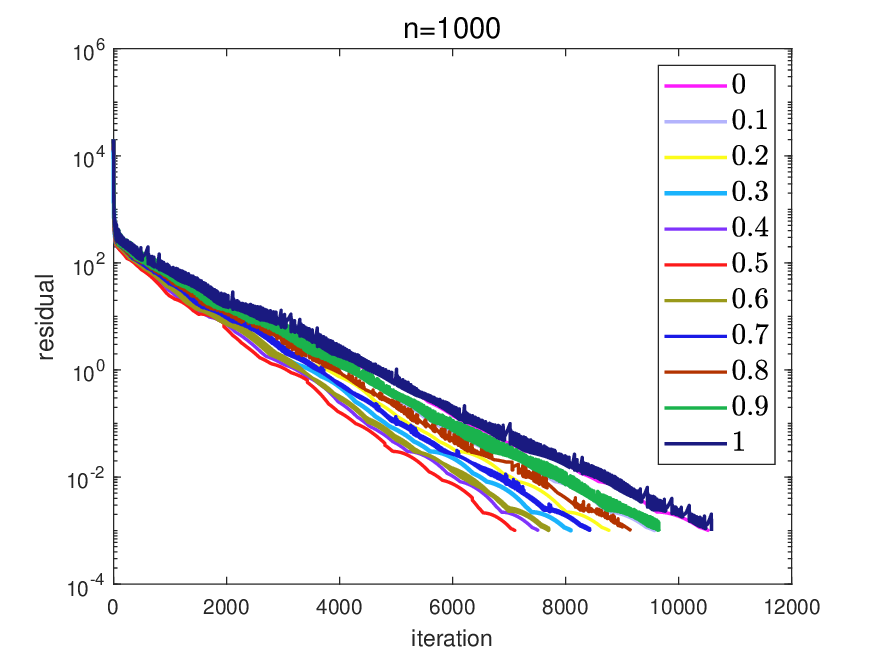}
	}
	\subfigure{
		\includegraphics[width=3.1in]{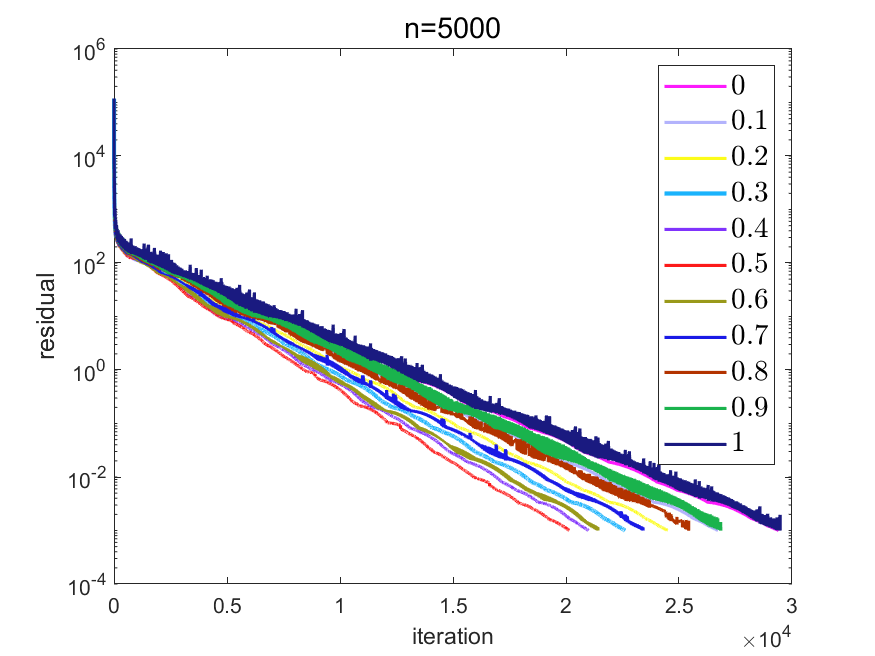}
	}
	\caption{Numerical performance of Algorithm 3 with different $\beta $ for Experiment 2}
		\label{fig2}
			\end{center}

\end{figure}

The comparisons of numerical results of Algorithm 3 under different dimensions and $\beta$ are shown in Table \ref{table2} and Figure \ref{fig2}. We note that when $\beta=0$, Algorithm 3 reduces to the gradient descent algorithm, and when $\beta=1$, Algorithm 3 reduces to the extra-gradient algorithm. It can be seen from these tables and figures that at various dimensions, the numerical performance is best when $\beta=0.5$, and improves by approximately 33\% compared to the gradient descent algorithm. The numerical performance for other values of $\beta$ is symmetrically distributed. This is consistent with our previous explanation from the perspective of differential equations, and demonstrates the advantage of incorporating a correction coefficient in the descent direction.

\vspace {-0.5cm}
\begin{center}
	\begin{table}[htbp]
		\centering
		\caption{Comparison of the efficiency of Algorithm 3 under different $\beta$ for Experiment 2}
		\begin{tabular}
			{m{1.2cm}<{\centering}|m{1.2cm}<{\centering}m{1.2cm}<{\centering}|m{1.2cm}<{\centering}m{1.2cm}<{\centering}|m{1.2cm}<{\centering}m{1.2cm}<{\centering}}
			\toprule
			&\multicolumn{2}{c|}{$n=1000$}&\multicolumn{2}{c|}{$n=2000$}&\multicolumn{2}{c}{$n=5000$}\\\hline
			&Iter.&CPU&Iter.&CPU&Iter.&CPU\\
			$\beta=0$&10534&2.3537&16004&21.2347&29406&325.5482\\
			$\beta=0.1$&9575&2.0596&14549&20.3756&26732&293.5786\\
			$\beta=0.2$&8776&1.8370&13336&17.7952&24504&271.1818\\
			$\beta=0.3$&8100&1.6974&12310&16.3450&22618&257.8097\\
			$\beta=0.4$&7520&1.6083&11429&15.2138&21002&233.4083\\
			$\beta=0.5$&\textbf{7106}&\textbf{1.5463}&\textbf{10750}&\textbf{14.5367}&\textbf{20149}&\textbf{224.5374}\\
			$\beta=0.6$&7704&1.8194&11688&16.9927&21414&254.4506\\
			$\beta=0.7$&8430&1.9288&12768&19.5998&23414&267.9760\\
			$\beta=0.8$&9159&2.0952&14011&19.3668&25432&301.5660\\
			$\beta=0.9$&9654&2.2999&14642&21.3543&26850&318.6725\\
			$\beta=1$&10583&2.4761&16037&22.9588&29497&339.8504\\
			\bottomrule
		\end{tabular}
		\label{table2}
	\end{table}
\end{center}

\section{Conclusion}
This paper explores prediction-correction algorithms for solving pseudo-convex and convex optimization problems. We express the descent direction of some typical prediction-correction algorithms as an adjustment to the gradient direction, and three improved prediction-correction algorithms are proposed by selecting suitable adjustment coefficients. These coefficients are intuitively explained as rotations of the descent direction, demonstrating that rotating by a suitable angle results in a more effective descent direction. Several specific adjustment coefficients are analyzed from the perspective of discrete differential equations. The ranges of these coefficients are provided for both pseudo-convex and convex objective functions, accompanied by a convergence analysis. The final numerical results indicate that the algorithm achieves optimal performance when the adjustment coefficient approximates the trapezoidal formula, which possesses higher-order global discrete error, underscoring the value and significance of this paper.

\bibliographystyle{unsrt}

\bibliography{cas-refs.bib}

@ARTICLE{ZJ15,
  author  = {Zappone A. and Jorswieck E.},
  title   = {Energy efficiency in wireless networks via fractional programming theory},
  journal = {Foundations Trends Commun. Inf. Theory.}, 
  volume  = {11},
  year    = {2015},
  pages   = {185-396}
}

@ARTICLE{AC20,
  author  = {Attouch H. Csaba L.},
  title   = {Newton-like inertial dynamics and proximal algorithms governed by maximally monotone operators},
  journal = {SIAM J. Optim.}, 
  volume  = {30},
  year    = {2020},
  pages   = {3252-3283}
}

@CONFERENCE{Be17,
	author  = {Bertsekas D.},
	title   = {Nonlinear programming. 3rd Edition},
	booktitle = {Athena Scientific, Belmont, Massachusetts}, 
	year    = {2017}
}

@CONFERENCE{BM19,
	author  = {Bo\c{t} R.I. Mertikopoulos P. Staudigl M. Vuong P.T.},
	title   = {On the convergence of stochastic forward-backward-forward algorithms with variance reduction in pseudo-monotone variational inequalities},
	booktitle = {NIPS 2018-Workshop on Smooth Games, Optimization and Machine Learning.}, 
	year    = {2018},
	pages   = {1-5}
}

@ARTICLE{CH13,
	author  = {Cai X. Han D. Xu L.},
	title   = {An improved first-order primal-dual algorithm with a new correction step},
	journal = {J. Global Optim.},
	volume  = {57}, 
	year    = {2013},
	pages   = {1419-1428}
}

@ARTICLE{Ca19,
	author  = {Cai X.},
	title   = {A proximal point algorithm with asymmetric linear term},
	journal = {Optim. Lett.},
	volume  = {13}, 
	year    = {2019},
	pages   = {777-793}
}

@ARTICLE{CG22,
	author  = {Cai X. Guo K. Jiang F.},
	title   = {The developments of proximal point algorithms},
	journal = {J. Oper. Res. Soc. China},
	volume  = {10}, 
	year    = {2022},
	pages   = {197-239}
}

@ARTICLE{CL17,
	author  = {Chen Y. Lan G. Ouyang Y.},
	title   = {Accelerated schemes for a class of variational inequalities},
	journal = {Math. Program.},
	volume  = {165}, 
	year    = {2017},
	pages   = {113-149}
}

@misc{CL21,
	author  = {Chen L. Luo H.},
	title   = {A unified convergence analysis of first order convex optimization methods via strong Lyapunov functions},
	howpublished = "\url{arXiv preprint arXiv:2108.00132v1}",
	year    = {2021}
}

@ARTICLE{CY13,
	author  = {Cheung K. Yang S. Hanzo L.},
	title   = {Achieving maximum energy efficiency in multi-relay OFDMA cellular networks: A fractional programming approach},
	journal = {IEEE Trans. Commun.},
	volume  = {61}, 
	year    = {2013},
	pages   = {2746-2757}
}

@ARTICLE{DO17,
	author  = {Diakonikolas J. Orecchia L.},
	title   = {The approximate duality gap technique: a unified theory of first-order methods},
	journal = {SIAM J. Optim.},
	volume  = {29}, 
	year    = {2019},
	pages   = {660-689}
}

@ARTICLE{DC18,
	author  = {Dong X. Cai X. Han D.},
	title   = {Prediction-correction method with {BB} step sizes},
	journal = {Front. Math. China},
	volume  = {13}, 
	year    = {2018},
	pages   = {1325-1340}
}

@CONFERENCE{FP03,
	author  = {Facchinei F. Pang J.},
	title   = {Finite-dimensional variational inequalities and complementarity problems. Vol. I},
    booktitle = {Springer, New York},
	year    = {2003}
}

@ARTICLE{JC21,
	author  = {Jiang F. Cai X. Han D.},
	title   = {The indefinite proximal point algorithms for maximal monotone operators},
	journal = {Optimization},
	volume  = {70}, 
	year    = {2021},
	pages   = {1759-1790}
}

@ARTICLE{HF09,
	author  = {He B. Fu X. Jiang Z.},
	title   = {Proximal-point algorithm using a linear proximal term},
	journal = {J. Optim. Theory Appl.},
	volume  = {141}, 
	year    = {2009},
	pages   = {299-319}
}

@ARTICLE{HL02,
	author  = {He B. Liao L.},
	title   = {Improvements of some projection methods for monotone nonlinear variational inequalities},
	journal = {J. Optim. Theory Appl.},
	volume  = {112}, 
	year    = {2002},
	pages   = {111-128}
}

@ARTICLE{HY04,
	author  = {He B. Yuan X. Zhang J.},
	title   = {Comparison of two kinds of prediction-correction methods for monotone variational inequalities},
	journal = {Comput. Optim. Appl.},
	volume  = {27}, 
	year    = {2004},
	pages   = {247-267}
}

@ARTICLE{HY13,
	author  = {He B. Yuan X. Zhang W.},
	title   = {A customized proximal point algorithm for convex minimization with linear constraints},
	journal = {Comput. Optim. Appl.},
	volume  = {56}, 
	year    = {2013},
	pages   = {559-572}
}

@CONFERENCE{JL18,
	author  = {Jelena D. Lorenzo O.},
	title   = {Accelerated extra-gradient descent: a novel accelerated first-order method},
	booktitle = {9th Innovations in Theoretical Computer Science Conference.},
	year    = {2018},
	pages   = {1-19}
}

@ARTICLE{Ko76,
	author  = {Korpelevi\v{c} G.},
	title   = {An extragradient method for finding saddle points and for other problems},
	journal = {Matecon.},
	volume  = {12}, 
	year    = {1976},
	pages   = {747-756}
}

@CONFERENCE{Lin05,
	author  = {Lin C.},
	title   = {Numerical computation methods (in Chinese)},
    booktitle = {Science Press.},
	year    = {2005}
}

@CONFERENCE{LH20,
	author  = {Liu H. Hu J. Li Y. Wen Z.},
	title   = {Optimization: modeling, algorithm and theory (in Chinese)},
	booktitle = {Higher Education Press.},
	year    = {2020}
}

@ARTICLE{MT20,
	author  = {Malitsky Y. Tam M.},
	title   = {A forward-backward splitting method for monotone inclusions without cocoercivity},
	journal = {SIAM J. Optim.},
	volume  = {30}, 
	year    = {2020},
	pages   = {1451-1472}
}

@CONFERENCE{Ma70,
	author  = {Mangasarian O.L.},
	title   = {Convexity, pseudo-convexity and quasi-convexity of composite functions},
	booktitle = {Stochastic Optimization Models in Finance, Academic Press},
	year    = {1975}
}

@ARTICLE{Ne05,
	author  = {Nemirovski A.},
	title   = {Prox-method with rate of convergence O(1/t) for variational inequality with
	Lipschitz continuous monotone operators and smooth convex-concave saddle point problems},
	journal = {SIAM J. Optim.},
	volume  = {15}, 
	year    = {2005},
	pages   = {229-251}
}

@CONFERENCE{Ne18,
	author  = {Nesterov Y.},
	title   = {Lectures on convex optimization},
	booktitle = {Springer Optimization and Its Applications.},
	year    = {2018}
}

@ARTICLE{OA22,
	author  = {Ogwo G.N. Alakoya T.O. Mewomo O.T.},
	title   = {Inertial iterative method with self-adaptive step size for finite family of split monotone 	variational inclusion and fixed point problems in Banach spaces},
	journal = {Demonstr. Math.},
	volume  = {55}, 
	year    = {2022},
	pages   = {193-216}
}

@ARTICLE{PB13,
author  = {Parikh N. Boyd S.},
title   = {Proximal algorithms},
journal = {Foundations and Trends$^{\circledR}$ in Optimization.},
volume  = {1}, 
year    = {2013},
pages   = {123-231}
}

@ARTICLE{SB16,
author  = {Su W. Boyd S. Cand\`{e}s E.},
title   = {A differential equation for modeling Nesterov's accelerated gradient method: theory and insights},
journal = {J. Mach. Learn. Res.},
volume  = {17}, 
year    = {2016},
pages   = {1-43}
}

@CONFERENCE{SM92,
author  = {Stancu-Minasian I. M.},
title   = {Fractional programming: theory, methods and applications},
booktitle = {Norwell, MA, USA: Kluwer.},
year    = {1992}
}

@ARTICLE{SY18,
author  = {Shen K. Yu W.},
title   = {Fractional programming for communication systems-part {I}: power control and beamforming},
journal = {IEEE Trans. Signal Process.},
volume  = {66}, 
year    = {2018},
pages   = {2616-2630}
}

@ARTICLE{Ts00,
author  = {Tseng P.},
title   = {A modified forward-backward splitting method for maximal monotone mappings},
journal = {SIAM J. Control Optim.},
volume  = {38}, 
year    = {2000},
pages   = {431-446}
}

@misc{ZB20,
author = {Zhang G. Bao X. Lessard L.},
title  = {A unified analysis of first-order methods for smooth games via integral quadratic constraints},
howpublished = "\url{arXiv preprint arXiv:2009.11359v4}",
year = {2021}
}

\end{document}